\documentclass[12pt]{article} 
\usepackage{amsmath, amscd, amssymb,latexsym, epsfig, color}
\input xy
\xyoption{all}

\newtheorem{theorem}{Theorem}[section]
\newtheorem{proposition}[theorem]{Proposition}
\newtheorem{question}[theorem]{Question}
\newtheorem{corollary}[theorem]{Corollary}

\newtheorem{conjecture}[theorem]{Conjecture}
\newtheorem{remark}[theorem]{Remark} 
\newtheorem{lemma}[theorem]{Lemma}

\newtheorem{problem}[theorem]{Problem}

\def\proof{\smallskip\noindent {\it Proof: \ }} 
\def\endproof{\hfill$\square$\medskip}
 \def\Z{\mathbb{Z}} 
 
\def\Q{\mathbb{Q}}
\def\R{\mathbb{R}}
\def\C{\mathbb{C}}
\def\m{{\mathfrak{M}}} 
\def\x{{\bf{x}}}

\def\S{\mathcal{S}} 
\def\H{\mathcal{H}}
\def\M{\mathcal{M}}

\newcommand{\inc}{\iota}
\newcommand{\con}{\phi}
\newcommand{\pr}{\pi}
\newcommand{\lk}{\mbox{\upshape lk}\,} 
 
\newcommand{\field}{{\bf k}}

\newcommand{\Init}{\mbox{\upshape In}\,} 
\newcommand{\Star}{\mbox{\upshape st}\,} 
\newcommand{\rk}{\mbox{\upshape rk}\,} 
 
\newcommand{\soc}{\mbox{\upshape Soc}\,} 
\newcommand{\Span}{\mbox{\upshape Span}\,}
\newcommand{\GL}{\mbox{\upshape GL}}
\newcommand{\Char}{\mbox{\upshape char}\,}

\title{Socles of Buchsbaum modules, complexes and posets.}
\author{Isabella Novik 
\thanks{Research partially supported by Alfred P.~Sloan Research
Fellowship and NSF grant DMS-0500748}\\ 
\small Department of Mathematics, Box 354350\\[-0.8ex] 
\small University of Washington, Seattle, WA 98195-4350, USA,\\[-0.8ex] 
\small \texttt{novik@math.washington.edu} 
\and Ed Swartz 
\thanks{Research partially supported by NSF grant DMS-0600502}\\ 
\small Department of Mathematics, \\[-0.8ex] 
\small Cornell University, Ithaca NY, 14853-4201, USA, \\[-0.8ex] 
\small \texttt{ebs22@cornell.edu } }

\begin{document}
\maketitle 

\begin{abstract} 

The socle of a graded Buchsbaum module is studied and is related
to its local cohomology modules. This algebraic result is then applied
to face enumeration of Buchsbaum simplicial complexes and posets.
In particular, new necessary conditions on face numbers and Betti numbers
of such complexes and posets are established. These conditions are
used to settle in the affirmative
K\"uhnel's conjecture for the maximum value of the Euler characteristic of 
a $2k$-dimensional simplicial manifold on $n$ vertices as well as Kalai's 
conjecture providing a lower bound on the number of edges of a simplicial 
manifold in terms of its dimension, number of vertices, and the first Betti 
number.
\end{abstract}

\section{Introduction}
This paper is devoted to the study of socles of graded Buchsbaum
modules and their applications to face enumeration.
A basic invariant of a simplicial complex or poset $\Delta$ 
is its $f$-vector $f(\Delta)=(f_0, f_1, \ldots, f_{d-1}),$
where $d-1$ is the dimension of $\Delta$ and $f_i$ is the number of 
its $i$-dimensional faces. One of the fundamental problems in 
geometric combinatorics is to characterize, or at least to 
obtain significant new
necessary conditions, on the $f$-vectors of various classes of complexes.
Here we study this problem for the class of Buchsbaum 
simplicial complexes and posets, and especially its subclass of 
complexes and posets representing manifolds.
We start by discussing the history of the problem and describing
our main results. All definitions are deferred to later sections.

Thirty years ago the pioneering work of Stanley and Hochster (see Chapter 2 
of \cite{St96}) made the study of combinatorics of simplicial complexes 
inseparable from the study of monomial ideals and graded algebras. 
Their insight was to associate with every simplicial complex a certain 
graded ring known today as the face ring or the Stanley-Reisner ring, 
and to read various combinatorial and topological invariants/properties 
of a complex off of the algebraic invariants of that ring. 

Call a simplicial complex $\Delta$ Cohen-Macaulay, resp.~Buchsbaum, 
if its Stanley-Reisner ring is Cohen-Macaulay, resp.~Buchsbaum.
Reisner \cite{Reisner}, building on (then unpublished) work of Hochster, 
gave a purely combinatorial-topological characterization of Cohen-Macaulay 
complexes, while Stanley worked out a complete characterization 
of $f$-vectors of Cohen-Macaulay complexes \cite{St77}, and later of 
$f$-vectors of Cohen-Macaulay simplicial posets \cite{St91}. 
In \cite{Sch}, Schenzel analyzed general Buchsbaum rings and modules and used
these algebraic results to generalize both Reisner's result
to a combinatorial-topological 
characterization of Buchsbaum complexes and the necessity portion
of Stanley's result to certain necessary conditions on the $f$-vectors and 
Betti numbers of Buchsbaum complexes.

One motivation for the study of
$f$-vectors of Cohen-Macaulay, resp.~Buchsbaum complexes
came from the desire to extend McMullen's upper bound theorem \cite{McMullen}
(UBT, for short) that provided sharp upper bounds on the face numbers of 
polytopes in terms of their dimension and the number of vertices 
to the class of simplicial spheres and, more generally,
Eulerian simplicial manifolds.
That such an extension does hold was conjectured by Klee \cite{Klee64-UBT}, 
and proved by Stanley \cite{St75} for the case of spheres, and then by 
Novik \cite{N98} for several classes of simplicial manifolds including all 
Eulerian ones. Novik's proof relied on Schenzel's results and on the method 
of algebraic shifting introduced by Kalai, see e.g.~\cite{Kalai02}. 
The main ingredient of the proof was a certain strengthening of Schenzel's 
conditions on the $f$-vectors and Betti numbers of Buchsbaum complexes.

In this paper we strengthen 
these conditions even further -- see 
Theorems \ref{main2} and \ref{new-bounds}, verifying 
in the affirmative a part of Kalai's conjecture, see
\cite[Conjecture 36]{Kalai02}. 
To derive these conditions we establish 
a new commutative algebra result, see
Theorems \ref{main-thm} and \ref{main-thm2}, 
that relates the socle of
a general Buchsbaum module to its local cohomology modules,
a result that we hope will be of interest in its own right. 
The same algebraic theorem is then used to
show that the lower bound part of our conditions 
on the $f$-vectors and Betti numbers
applies also to all Buchsbaum simplicial posets, 
see Theorem \ref{main-posets}. Based on the situation in 
dimensions up to four (see Section 7), we believe that these lower bound
conditions provide, in fact, a complete characterization of the $f$-vectors of
Buchsbaum simplicial posets with prescribed Betti numbers.

Related to the UBT is a conjecture 
by  K\"uhnel \cite[Conjecture B]{Kuh95}
for the maximum value of the Euler characteristic of a 
$2k$-dimensional simplicial
manifold on $n$ vertices. This conjecture was previously known to hold 
only for manifolds
with at least $4k+3$ or  at most $3k+3$ vertices \cite{N98, N05}. 
Here we use  Theorem \ref{new-bounds} to
prove it for all values of $n$, see Theorem \ref{Kuhnel-conj}.

In \cite{Kalai87}, Kalai conjectured a lower bound on the number of edges of 
a simplicial manifold in terms of its dimension, number of vertices, 
and the first Betti number. This conjecture was verified by Swartz \cite{Sw}
for manifolds whose first Betti number is one 
as well as for orientable manifolds of dimension at least four
with  vanishing second Betti number.  
Here we prove this conjecture for all  manifolds, 
see Theorem \ref{g2}.

The structure of the paper is as follows. 
In Section 2 after providing a necessary 
background on Buchsbaum modules, we state and prove our main algebraic result, 
Theorem \ref{main-thm}, on which
all other theorems of this paper are based. Section 3 contains an overview of 
simplicial complexes and their Stanley-Reisner rings as well 
as a combinatorial-topological
translation of Theorem \ref{main-thm} for the case of 
Buchsbaum simplicial complexes. 
Section 4 is devoted to deriving new upper bounds on the face numbers of 
Buchsbaum simplicial complexes and in particular includes 
the proof of the K\"uhnel 
conjecture. In Section 5 we prove Kalai's lower bound conjecture. 
In Sections 6 we study $f$-vectors of Buchsbaum simplicial posets. 
In Section 7, we discuss several examples and state 
a few  of the many still open questions on the $f$-vectors of 
Buchsbaum simplicial complexes and posets.

\section{Socles in terms of local cohomology} 
In this section we state and prove our main algebraic result concerning the 
socle of a Buchsbaum module. This theorem is the key to all the combinatorial
applications discussed in the rest of the paper.
 
We start by reviewing necessary background material. 
For all undefined terminology 
as well as for more details the reader is referred to \cite{StVo}.
Let $\field$ be an {\bf infinite} field of an arbitrary characteristic 
and let $S:=\field[x_1, \ldots, x_n]$ be a polynomial ring. 
We denote by $\m$ the irrelevant ideal of $S$, and by $\m_j$ the 
$j$th homogeneous component of $\m$.
All modules considered in this paper are Noetherian ($\Z$-)graded 
modules over $S$.

Let $M$ be a Noetherian graded $S$-module of Krull dimension $d\geq 0$. 
A {\em homogeneous system of parameters} of $M$, abbreviated h.s.o.p,  
is a sequence $\theta_1, \cdots, \theta_d$ of homogeneous elements of $\m$
such that $\dim M/(\theta_1, \cdots, \theta_d)M=0$  
(equivalently, $M/(\theta_1, \cdots, \theta_d)M$ 
is a finite-dimensional vector space over \field).
A h.s.o.p.~all of whose elements belong to $\m_1$ is called a 
{\em linear system of parameters}, 
l.s.o.p. for short. It follows from the Noether Normalization Lemma that a 
l.s.o.p. always exists.
A sequence of elements $\theta_1, \cdots, \theta_r \in \m$ is a 
{\em weak $M$-sequence} if for each 
$i=1, \ldots, r$
$$
(\theta_1, \cdots, \theta_{i-1})M : \theta_i = 
(\theta_1, \cdots, \theta_{i-1})M : \m.
$$

Our main object of study is the class of {\em Buchsbaum modules}. 
Following Definition~3.1 on page 95 of \cite{StVo} combined with Theorem 3.7 
on page 97,
we say that a Noetherian graded $S$-module $M$ is {\em Buchsbaum} 
if every h.s.o.p.~of $M$ is a weak $M$-sequence. 
Since any regular $M$-sequence is also a weak $M$-sequence, 
all Cohen-Macaulay modules are Buchsbaum. 
A large family of Buchsbaum modules most of which are not Cohen-Macaulay 
is given by the face rings 
of triangulated manifolds --- see Section \ref{St-R-rings}.

The following lemma summarizes several basic properties 
of Buchsbaum modules we will rely on frequently. 
Here $H^i(M)$ denotes the $i$th local cohomology of $M$ with respect to $\m$. 
In particular, 
$$H^0(M)=0 : \m^\infty =\{y \in M \,|\, \m^k y =0 \mbox{ for some } k>0\}$$ 
is a submodule of $M$. The modules $H^i(M)$ are graded provided $M$ is.

\begin{lemma} \label{Buchs-modules-properties}
Let $M$ be a graded Noetherian $S$-module of Krull dimension $d\geq 0$.
If $M$ is Buchsbaum and $\theta_1, \ldots, \theta_r$ is a part of a 
h.s.o.p.~for $M$, then
\begin{enumerate}
\item $M/(\theta_1, \ldots, \theta_r)M$ is a Buchsbaum module of Krull 
dimension $d-r$,
\item $(\theta_1, \ldots, \theta_{r-1})M : \m=
(\theta_1, \ldots, \theta_{r-1})M : \theta_r = 
(\theta_1, \ldots, \theta_{r-1})M : \theta_r^2$, and
\item $\m \cdot H^i(M)=0$ for all $0\leq i<d$.
\end{enumerate}
\end{lemma}
\noindent All parts of the lemma can be found in \cite{StVo}: 
for (1) see Corollary 1.11 on page 65, 
for (2) see Proposition 1.10 on pages 64-65, and for (3)
 see Corollary 2.4 on page 75.

Recall that the {\em socle} of a module $M$ is 
$$ \soc M:= 0 : \m = \{y\in M \,|\, \m \cdot y=0\}.$$
We are now in a position to state our main theorem relating the socle  
to the local cohomology modules. 
We denote by $M_k$ the $k$th homogeneous component of a graded module $M$, 
and by $rM$ the direct sum of $r$ copies of $M$. 

\begin{theorem} \label{main-thm}
Let $M$ be a Noetherian graded $S$-module of Krull dimension $d$, and let
$\theta_1, \ldots, \theta_d$ be a l.s.o.p. If $M$ is Buchsbaum, 
then for all integers $k$,
$$
\left(\soc M/(\theta_1, \ldots, \theta_d)M \right)_k \cong 
\left(\bigoplus_{j=0}^{d-1} {d \choose j} H^j(M)_{k-j} \right) 
\bigoplus \S_{k-d},
$$ where $\S$ is a  graded submodule of $\soc H^d(M)$.
\end{theorem}

We begin the proof with the following lemma.  
For a graded module $M$, $M(-a)$ is the same module, 
but with grading $M(-a)_k = M_{k-a}.$

\begin{lemma}  \label{zero-map}
If $M$ is a Buchsbaum $S$-module and $\theta$ is a part of a l.s.o.p.~for $M$, 
then 
\begin{enumerate}
\item $H^i(\theta M)\cong H^i(M(-1))$ for all $i>0$, and 
\item the map 
$\inc^\ast: H^i(\theta M) \rightarrow H^i(M)$ induced by inclusion 
$\inc: \theta M \hookrightarrow M$ is the zero map 
for all $0\leq i < \dim M$.
\end{enumerate}
\end{lemma}

\proof We verify both claims simultaneously.
Let $N=M/H^0(M)$. Consider the following diagram and the induced local 
cohomology diagram:
$$
\begin{CD}
\theta M @>\inc>> M \\
         @VfVV      @VVqV \\
 N(-1) @>\cdot\theta>> N
\end{CD} 
\qquad \qquad \qquad
\begin{CD}
H^i(\theta M) @>\inc^{\ast}>> H^i(M) \\
              @Vf^{\ast}VV      @VV{q^\ast}V \\
H^i(N(-1)) @>\cdot\theta>> H^i(N).
\end{CD}  
$$
Here $\inc$ is the inclusion map, 
$q: x \mapsto x+H^0(M)$ is the quotient map, and 
$f$ is the map defined by $\theta x \mapsto x+ H^0(M)$. 
To see that $f$ is well-defined, suppose that 
$\theta x =\theta y$. Then $\theta (x-y)=0$, and
so, by the definition of a Buchsbaum module, 
$\m \cdot (x-y)=0$. Hence $x-y \in H^0(M)$.
We claim that $f$ (and hence also $f^\ast$) is an isomorphism. 
It is  surjective since for any $x \notin H^0 (M)$, $\theta x \neq 0.$ 
To show injectivity, assume
$x\in H^0(M)$. Then $\m^l\cdot x =0$ for some $l>0$. 
In particular, $\theta^{l} x=0$, which, by 
Part (2) of Lemma~\ref{Buchs-modules-properties}, 
implies that $\theta x =0$.  

The map $f$ was chosen to make the diagram on the left commute. 
By naturality of  local cohomology, 
the induced diagram also commutes. Now, since $H^0(M)$ has Krull dimension 0, 
$H^i(H^0(M))=0$ for all $i>0$,
and so $q^{\ast}$ is an isomorphism for $i>0$. If also $i<\dim M$, 
then Part (3) of Lemma \ref{Buchs-modules-properties} implies that 
the bottom horizontal map in the induced diagram is the zero map, and
we conclude that $\inc^\ast=0$ for all $0<i<\dim M$.
For $i=0$, another application of Part (2) of 
Lemma \ref{Buchs-modules-properties} 
shows that $H^0(\theta M) = 0$, and so $\inc^\ast=0$ in this case as well. 
This completes the proof of the second claim, while the string of isomorphisms
$$H^i(\theta M) \stackrel{f^\ast}{\longrightarrow}
 H^i(N(-1)) \stackrel{(q^\ast)^{-1}}{\longrightarrow} H^i(M(-1)) 
\quad \mbox{ for } i>0$$
implies the first claim.
\endproof

It is convenient to introduce the following notation:
for subsets $A, C$ of $\{1,2,\ldots, d\}=[d]$, 
write $\theta^C$ to denote $\prod_{i\in C} \theta_i$, and write $M(A)$ 
to denote 
$M/(\theta_i : i \notin A)M$. 
In particular, $M(\emptyset)=M/(\theta_1, \ldots, \theta_d)M$ and $M[d]=M$.
By repeated application of Lemma \ref{zero-map}(1), 
to prove Theorem \ref{main-thm} it 
is enough to verify the following. 
(We distinguish between strict and non-strict inclusions
by using symbols `$\subset$' and `$\subseteq$', respectively.)

\begin{theorem} \label{main-thm2}
  Let $M$ be a Noetherian graded $S$-module  of Krull dimension $d$, and let
$\theta_1, \ldots, \theta_d$ be a l.s.o.p. If $M$ is Buchsbaum, then 
$$
\soc(M(\emptyset)) \cong 
\left(\bigoplus_{C\subset [d]} H^{|C|}(\theta^C M) \right) 
\bigoplus \S,
$$ where $\S$ is a  graded submodule of $\soc H^d(\theta^{[d]}M)$.
\end{theorem}

The proof of Theorem \ref{main-thm2} involves ``chasing" 
a few commutative diagrams.
Our starting point is the short exact sequence
$$
0 \rightarrow \theta_s\cdot\theta^{C}M(A) 
\stackrel{\inc}{\rightarrow} \theta^{C}M(A) 
\stackrel{\pr_s}{\rightarrow} \theta^{C} M(A\setminus s)\rightarrow 0,
$$
where $A \subseteq [d], C \subset A$, $s\in A\setminus C$, 
$\inc$ is the inclusion map, and $\pr_s$ is the projection map. 
(The subscript $s$ indicates that $\pr_s$ maps a module to its 
quotient by the submodule generated by $\theta_s$.)
We refer to such a sequence as an $(A,C, s)$-sequence.
It gives rise to the long exact local cohomology sequence, 
where we denote by $\con_s^{\ast}$ the connecting homomorphism:
$$\cdots \rightarrow H^i(\theta^{C\cup s}M(A))  
\stackrel{\inc^{\ast}}{\rightarrow} H^i(\theta^{C}M(A))
\stackrel{\pr_s^{\ast}}{\rightarrow} H^i(\theta^{C}M(A\setminus s)) 
\stackrel{\con_s^{\ast}}{\rightarrow} H^{i+1}(\theta^{C\cup s}M(A))
\stackrel{\inc^{\ast}}{\rightarrow} \cdots.
$$

If $M$ is Buchsbaum, then, as follows from 
Lemma~\ref{Buchs-modules-properties}(1),
 $M(A)$ is also Buchsbaum and has Krull dimension $|A|$. 
Thus, by Lemma \ref{zero-map}, 
$\inc^{\ast}: H^i(\theta^{C\cup s}M(A))   
{\rightarrow} H^i(\theta^{C}M(A))$ is the zero map provided 
$i < |A|$. The above long exact sequence then breaks into the 
following short exact sequences:
\begin{eqnarray}  \label{short-seq}
 0 \rightarrow  H^i(\theta^{C}M(A )) 
\stackrel{\pr_s^{\ast}}{\hookrightarrow} H^i(\theta^{C}M(A \setminus s)) 
\stackrel{\con_s^{\ast}}{\rightarrow} H^{i+1}(\theta^{C\cup s}M(A ))
\rightarrow 0,  \;\; \mbox{ for all }i<|A|-1 \\
 \label{last-seq}
 0 \rightarrow  H^{|A|-1}(\theta^{C}M(A))
\stackrel{\pr_s^{\ast}}{\hookrightarrow} H^{|A|-1}(\theta^{C}M(A \setminus s)) 
\stackrel{\con_s^{\ast}}{\rightarrow} H^{|A|}(\theta^{C\cup s}M(A)) 
\stackrel{\inc^{\ast}}{\rightarrow}  H^{|A|}(\theta^{C}M(A)). 
\end{eqnarray}

For $A = \{1\}$, $s=1$, and  $C=\emptyset,$ (\ref{last-seq}) becomes
\begin{equation} \label{A=d}
0 \rightarrow  H^0( M(\{1\}))
\stackrel{\pr_1^{\ast}}{\hookrightarrow} H^0( M(\emptyset)) 
\stackrel{\con_1^{\ast}}{\rightarrow} H^1(\theta_1 M(\{1\})) 
\stackrel{\inc^{\ast}}{\rightarrow}  H^1( M(\{1\})).
\end{equation}
Since  the Krull dimension of $M(\{1\})$ is one, the image of 
$\pr^\ast_1$ is contained in the socle of $H^0(M(\emptyset))
= M/(\theta_1,\dots,\theta_d) M.$  
This submodule of the socle can be analyzed using (\ref{short-seq}).

\begin{lemma}  \label{lhs}
  Let $M$ be a Buchsbaum module of Krull dimension $d$. 
If $C \subset A \subseteq [d]$, and $i < |A|,$ then
  \begin{equation} \label{easypart}
  H^i (\theta^C M(A)) \cong  \bigoplus_{D\subseteq [d]\setminus A} 
  H^{i+|D|}(\theta^{C \cup D} M).
  \end{equation}
\end{lemma}

\proof
The proof is by induction on $d-|A|$.  
If $d-|A|=0,$ then  $A=[d]$, and the lemma is equivalent to
$H^i (\theta^C M) \cong H^i (\theta^C M)$.
For the induction step, the short exact sequence~(\ref{short-seq}) 
implies that for $s\in[d]\setminus A$
$$
H^i (\theta^C M(A)) \cong H^i(\theta^C M(A \cup s)) 
\oplus H^{i+1}(\theta^{C \cup s} M(A \cup s)).
$$
The induction hypothesis applied to the two terms on 
the right-hand side finishes the proof. 
\endproof

\begin{corollary}  \label{lhs-socle-cor}
Let $M$ be a Buchsbaum module of Krull dimension $d$.  Then
$$ H^0(M(\{1\}) \cong  
\bigoplus_{C\subseteq [d]\setminus\{1\}} H^{|C|}(\theta^C M).$$ 
\end{corollary}

In view of (\ref{A=d}), we  have accounted for those terms 
of the direct sum in Theorem \ref{main-thm2} such that 
$1 \notin C.$   To finish the proof we examine 
the image of the socle of $H^{|A|-1}(M(A \setminus s))$ 
under $\con^\ast_s$ in (\ref{last-seq}), 
then specialize to the case $A=\{1\}$ and $s=1.$   

If $r\in A\setminus (C \cup s)$, then the $(A,C,s)$-sequence and the 
$(A\setminus r,C,s)$-sequence can be combined together to form 
the following commutative diagram:
$$\minCDarrowwidth22pt\begin{CD}
0 @>>>   \theta^{C\cup s}M(A\setminus r) 
 @>{\inc}>> \theta^{C}M(A\setminus r) 
 @>{\pr_s}>> \theta^{C} M(A\setminus\{r,s\}) @>>> 0 \\
@. @AA{\pr_r}A  @AA{\pr_r}A  @AA{\pr_r}A @.\\
0 @>>>   \theta^{C\cup s}M(A) 
 @>{\inc}>> \theta^{C}M(A) 
 @>{\pr_s}>> \theta^{C} M(A\setminus s) @>>> 0. \\
\end{CD}
$$
Naturality of local cohomology then implies that the diagram whose 
rows consist of 
the corresponding long exact local cohomology sequences 
and all of whose vertical maps are induced by $\pr_r$ also commutes. 
This observation together with equations (\ref{short-seq}) 
and (\ref{last-seq}) yields the following.

\begin{lemma} \label{commute}
Let $M$ be a Buchsbaum module of Krull dimension $d$. 
If $C \subset A \subseteq [d]$,
$s\in A\setminus C$, and  $r\in A\setminus ( C \cup s)$, 
then for all $i<|A|-1$ we have the following 
commutative diagram whose rows are exact
$$
\xymatrix{
H^{i}(\theta^{C}M(A\setminus r))\ar@{^{(}->}[r]^{\pr_s^{\ast}}  & H^ 
{i}(\theta^{C}M(A\setminus \{r,s\})) \ar[r]^{\con_s^{\ast}} & H^{i+1} 
(\theta^{C\cup s}M(A\setminus r)) & \\
&  H^{i}(\theta^{C} M(A\setminus s)) \ar@{^{(}->}[u]_{\pr_r^{\ast}}  
\ar@{->>}[r]^{\con_s^{\ast}} & H^{{i}+1}( \theta^{C\cup s}M(A))\ar[r]  
\ar@{^{(}->}[u]_{\pr_r^{\ast}} & 0
}
$$
\end{lemma}

Recall that our goal is to compute the image of 
$\soc H^{|A|-1}(M(A \setminus s))$ under $\con_s^\ast$ in (\ref{last-seq}).
We do this by induction with 
the following lemma allowing the inductive step.

\begin{lemma} \label{ind-step}
With the assumptions of Lemma \ref{commute}, 
for all $i<|A|-1$, the preimage of 
$\pr_r^{\ast} \left(H^{{i}+1}( \theta^{C\cup s}M(A))\right)$ 
under $\con_s^\ast$ is contained in the socle 
of $H^{i}(\theta^{C}M(A\setminus\{r,s\}))$.
\end{lemma}

\proof Let $y\in H^{{i}+1}(\theta^{C\cup s}M(A))$. 
By surjectivity of $\con_s^\ast$ (see the diagram of Lemma \ref{commute}),
there exists $x\in H^{i}(\theta^{C} M(A\setminus s))$ 
such that $\con_s^\ast(x)=y$. Since $\dim(M(A\setminus s))=|A|-1 > i$
it follows from Lemma \ref{Buchs-modules-properties}, 
that all elements of $H^{i}(\theta^{C} M(A\setminus s))$, 
including $x$, are in the socle of 
$H^{i}(\theta^{C} M(A\setminus s))$, and hence 
$\pr_r^\ast(x)\in \soc H^{i}(\theta^{C}M(A\setminus\{r,s\}))$.
But the diagram of Lemma \ref{commute} commutes, and so
$\pr_r^\ast(x) \in (\con_s^\ast)^{-1}(\pr_{r}^\ast (y))$. 
We have proved that each element $y\in H^{{i}+1}(\theta^{C\cup s}M(A))$ 
has a 
representative $\tilde{y}\in (\con_s^\ast)^{-1}(\pr_{r}^\ast (y))$ 
that lies in 
$\soc H^{i}(\theta^{C}M(A\setminus\{r,s\}))$. 
Now choose a \field-basis $B=\{y_1, \ldots, y_t\}$ for
$H^{{i}+1}(\theta^{C\cup s}M(A))$ and let 
$\tilde{B}=\{\tilde{y_1}, \ldots, \tilde{y_t}\}$ 
be a set of their representatives in the pull-back that lie in the socle. 
By Lemma \ref{commute},
$$(\con_s^\ast)^{-1} 
\pr_r^{\ast} \left(H^{{i}+1}( \theta^{C\cup s}M(A))\right)=
\pr_s^\ast(H^{i}(\theta^{C}M(A\setminus r))) \oplus \Span(\tilde{B})  ,$$
 and the assertion follows, since $\Span(\tilde{B})$ 
is in the socle by the choice of $\tilde{B}$,
and the first summand of the above decomposition 
is also in the socle by Lemma \ref{Buchs-modules-properties}(3).
\endproof

Using a fixed $r$ and varying $s$'s in Lemma \ref{commute}, 
we can chain the corresponding commutative squares together to obtain that
for subsets $T=T_k=\{s_1, \ldots, s_k\}$, $A$, and $C$ of $[d]$ satisfying 
$T\subseteq A\setminus C$ and $r\in A\setminus T$, and for $i<|A|-k$,
the following diagram with 
$\con_{T}^\ast:=\con_{s_1}^\ast \circ \cdots \circ \con_{s_k}^\ast$ commutes:

$$
\xymatrix{
\ker \con^\ast_T \ar@{^{(}->}[r]  & H^ 
{i}(\theta^{C}M(A\setminus (T \cup r)) \ar[r]^{\con_T^{\ast}} & H^{i+k} 
(\theta^{C\cup T}M(A\setminus r)) & \\
&  H^{i}(\theta^{C} M(A\setminus T)) \ar@{^{(}->}[u]_{\pr_r^{\ast}}  
\ar@{->>}[r]^{\con_T^{\ast}} & H^{{i}+k}( \theta^{C\cup T}M(A))\ar[r]  
\ar@{^{(}->}[u]_{\pr_r^{\ast}} & 0
}
$$
Here $\ker\con^\ast_{T}=
\oplus_{j=1}^k (\con_{s_{j+1}}^\ast \circ \ldots \circ \con^\ast_{s_k})^{-1} 
\pr_{s_j}^\ast
H^{i+k-j}(\theta^{C\cup\{s_{j+1}, \ldots, s_k\}} 
M(A\setminus(T_{j-1}\cup r)))$, so
the same argument as in Lemma \ref{ind-step} plus induction on $k$ then implies

\begin{lemma}  \label{ind-step-comb}
For all $i<|A|-|T|$, the preimage of 
$\pr_r^\ast H^{i+|T|}(\theta^{C\cup T}M(A))$ under 
$\con_T^\ast$ is contained in 
the socle of $H^i(\theta^C M(A\setminus(T\cup r)))$.
\end{lemma}

Now consider the following diagram.
$$\begin{array}{ccc}
 &  &H^0(M(\emptyset)) \\
          & &  \downarrow \con^\ast_1 \\
H^1(\theta^{[1]} M([2])) &\stackrel{\pr^\ast_{2}}{\hookrightarrow} 
            & H^1(\theta^{[1]}M([1]))\\
  & &    \downarrow \con^\ast_2 \\
  H^2(\theta^{[2]} M([3])) 
              &\stackrel{\pr^\ast_3}{\hookrightarrow} 
               & H^2(\theta^{[2]} M([2]))\\
& & \  \downarrow \con^\ast_3 \\
& &\vdots \\
& & \ \ \ \ \downarrow \con^\ast_{d-1}\\
H^{d-1} (\theta^{[d-1]} M([d])) 
              & \stackrel{\pr^\ast_d}{\hookrightarrow} 
          & H^{d-1} (\theta^{[d-1]} M([d-1]) \\
& & \ \ \downarrow \con^\ast_d\\
& & H^d(\theta^{[d]} M([d])).
\end{array}  $$
Lemma \ref{ind-step-comb} shows 
that for all $j,\  1 \le j \le d-1, $
$$(\con^\ast_{j} \circ \dots \circ \con^\ast_1)^{-1} 
\pr^\ast_{j+1} \ H^{j}(\theta^{[j]} M([j+1]))$$
 lies in the socle of $H^0(M(\emptyset)).$  
Using Lemma \ref{lhs} on each 
$H^{j}(\theta^{[j]} M([j+1]))$ 
accounts for all of the terms of the direct 
sum decomposition in Theorem \ref{main-thm2} 
with $1 \in C.$ 

Setting 
$\S = \con^\ast_d \circ \dots \circ \con^\ast_1 (\soc H^0(M(\emptyset))$ 
finishes the proof of Theorem \ref{main-thm2}.
 \endproof

\begin{remark}{\rm
Instead of graded Buchsbaum $S$-modules, one can work in the generality of 
Buchsbaum modules over Noetherian
local rings, see Definition 1.5 on page 63 
in \cite{StVo}.
A proof identical to that of Theorem \ref{main-thm2} 
then shows that if $M$ is a 
Noetherian module of Krull dimension $d$ over a local ring $A$, and 
$\theta_1, \ldots, \theta_d$ is a system of parameters of $M$, then
$$
\soc(M(\emptyset)) \cong 
\left(\bigoplus_{C\subset [d]} H^{|C|}(\theta^C M) \right) 
\bigoplus \S,
$$
provided $M$ is a Buchsbaum module.
Here $\S$ is a certain submodule of $\soc H^d(\theta^{[d]}M)$.
}
\end{remark}

\section{Simplicial complexes and Stanley-Reisner rings} \label{St-R-rings}
This section provides a short overview of several concepts 
and results related to simplicial complexes and their Stanley-Reisner rings.
A comprehensive reference to this topic is Chapter 2 of \cite{St96}.
The section concludes with a combinatorial-topological translation 
of Theorem~\ref{main-thm} for the case of Buchsbaum complexes and 
resulting new lower bounds on their face numbers.

Recall that a {\em simplicial complex} $\Delta$ on the vertex set $[n]$ is a 
collection of subsets of $[n]$ that is closed under inclusion and contains 
all singletons
$\{i\}$ for $i\in[n]$. The elements of $\Delta$ are called {\em faces}. 
The maximal faces (with respect to inclusion) are called {\em facets}. 
The {\em dimension of a face} $F\in \Delta$ is $\dim F:=|F|-1$. 
The {\em dimension of $\Delta$} is then defined as the maximal dimension 
of its faces.
A simplicial complex is called {\em pure} if all its facets have the same
 dimension.

We also need the notion of the link of a face: 
if $\Delta$ is a simplicial complex and $F$ is a face of $\Delta$,
then the {\em link} of $F$ in $\Delta$, $\lk(F)$, is the following subcomplex 
of $\Delta$
$$
\lk_\Delta(F):=\{G\in \Delta \, | \, G\cap F=\emptyset
 \mbox{ and } G\cup F \in \Delta\}.
$$  
Thus the link of the empty face is the complex itself.

A basic combinatorial invariant of a simplicial complex 
$\Delta$ on the vertex set $[n]$ is its {\em $f$-vector},
$f(\Delta)= (f_{-1}, f_0, \ldots, f_{d-1})$. Here $d-1=\dim \Delta$ 
and $f_i$ denotes
the number of $i$-dimensional faces of $\Delta$. In particular, 
$f_{-1}=1$ (there is only one empty face) and $f_0=n$.
An invariant that contains the same information as the $f$-vector, 
but sometimes is more convenient to work with, is the {\em $h$-vector}
 of $\Delta$,
$h(\Delta)=(h_0, h_1, \ldots, h_d)$ whose entries are defined by the 
following relation:
\begin{equation}   \label{h-vector}
\sum_{i=0}^d h_i x^{d-i} = \sum_{i=0}^d f_{i-1} (x-1)^{d-i}.
\end{equation}

The Stanley-Reisner ring of a simplicial complex
provides an important algebraic tool for studying $f$-numbers of 
simplicial complexes.
If $\Delta$ is a simplicial complex on $[n]$, then  
its {\em Stanley-Reisner ring} (also called the {\em face ring}) is 
$$\field[\Delta] = S/I_\Delta := \field[x_1, \ldots, x_n]/I_\Delta, 
     \quad \mbox{where }
I_\Delta=(x_{i_1}x_{i_2}\cdots x_{i_k} : \{i_1<i_2<\cdots<i_k\}\notin\Delta).
$$
The ideal $I_\Delta$ is called the {\em Stanley-Reisner ideal of $\Delta$}. 
(As in the previous section, here and throughout the paper, we assume that 
$\field$ is an infinite field of an arbitrary characteristic.) 
Since $I_\Delta$ is a monomial ideal, defining $\deg (x_i)=1$ 
for all $1\leq i \leq n$
makes $\field[\Delta]$ into a $\Z$-graded ring, 
while defining $\deg (x_i)=e_i$, 
where $e_1, \ldots, e_n$ is the standard basis for $\Z^n$, 
makes $\field[\Delta]$ into a $\Z^n$-graded ring.

The utmost significance of  Stanley-Reisner rings in the theory 
of $f$-numbers is explained by the fact that many 
combinatorial and topological properties of $\Delta$ translate to certain 
algebraic properties of $\field[\Delta]$ and vice versa. For instance, 
the Krull dimension of $\field[\Delta]$ 
(considered as a module over itself or over $S$) equals 
$\dim \Delta +1$, while 
the {\em ($\Z$-)Hilbert series of $\field[\Delta]$}, 
$F(\field[\Delta], x):= \sum_{j=0}^{\infty} \dim_\field \field[\Delta]_j x^j$ 
can be expressed in terms of the $h$-vector of $\Delta$:
\begin{equation} \label{Hilb-ser}
F(\field[\Delta], x) = (1-x)^{-d} \sum_{i=0}^d h_i x^i, 
\quad \mbox{where } d=\dim \Delta+1.
\end{equation}
(Both results can be found in \cite{St75} 
or on  pages 33, 54, and 58 of \cite{St96}.)
Moreover, the local cohomology of $\field[\Delta]$ 
(as a module over itself or over $S$) 
has a simple expression in terms of simplicial homology of
 the links of the faces of $\Delta$.
This result is known as  {\em Hochster's formula}, 
see \cite[Theorem II.4.1]{St96}. 

\begin{theorem} ({\bf Hochster})    \label{Hochster}
For a simplicial complex $\Delta$, $\alpha\in \Z^n$, 
$F=\{j\in [n] \, | \, \alpha_j\neq 0\}$,
and $i\geq 0$,
$$
H^i(\field[\Delta])_\alpha \cong \left\{\begin{array}{ll} 
         0, & \mbox{ if $F\notin\Delta$ or $\alpha_j>0$ for some $j\in[n]$}\\
         \widetilde{H}_{i-|F|-1} (\lk F; \field), & \mbox{ otherwise, }\\
\end{array}
\right.
$$
where $\widetilde{H}_{i}$ denotes the $i$th reduced simplicial homology 
with coefficients in $\field$.
\end{theorem}

Among the main objects of this paper are 
Cohen-Macaulay and Buchsbaum simplicial complexes.
A simplicial complex $\Delta$ is called {\em Cohen-Macaulay} (over \field), 
if $\field[\Delta]$ is Cohen-Macaulay 
(considered as a module over itself or over $S$).
Similarly, $\Delta$ is called {\em Buchsbaum} (over \field), 
if $\field[\Delta]$  is Buchsbaum.

Using Hochster's formula, Reisner \cite{Reisner} gave a purely 
combinatorial-topological characterization of Cohen-Macaulay complexes. 
His criterion was later generalized by Schenzel \cite{Sch} 
to the class of Buchsbaum complexes. We combine both these results in 
the following theorem.

\begin{theorem} \label{Reisner-Schenzel}
Let $\Delta$ be a simplicial $(d-1)$-dimensional complex. 
Then $\Delta$ is Cohen-Macaulay (over \field) if and only if 
$\widetilde{H}_{i} (\lk F; \field) = 0$ for all $F\in \Delta$, 
including $F=\emptyset$, and all
$i< d-|F|-1$. $\Delta$ is Buchsbaum (over \field) if and only if it 
is pure and the link of each vertex is Cohen-Macaulay (over \field).
\end{theorem}
Thus, every simplicial complex whose 
geometric realization is a \field-homology sphere 
is Cohen-Macaulay over \field, while a simplicial complex 
whose geometric realization is a \field-homology 
manifold is Buchsbaum over \field.
In the rest of the paper we refer to such complexes simply as 
homology spheres and homology manifolds, respectively.
Note that the class of homology spheres includes 
all triangulations of topological spheres. Similarly,
the class of homology manifolds includes all 
triangulations of topological manifolds. We refer to those complexes 
as simplicial spheres and simplicial manifolds, respectively.

In contrast to the $f$-numbers, the $h$-numbers may in general be negative.
The following result of Stanley, which is an immediate consequence of 
\cite[p.~35]{St96} and (\ref{Hilb-ser}), shows that this is not the case, 
if $\Delta$ is Cohen-Macaulay. For the rest of the paper we denote by
$(\Theta)$ the ideal of $\field[\Delta]$ generated by the elements
$\theta_1, \ldots, \theta_d$ of a l.s.o.p.~for $\field[\Delta]$.
\begin{theorem} ({\bf Stanley})   \label{CMh-numbers}
Let $\Delta$ be a $(d-1)$-dimensional Cohen-Macaulay complex 
and let $\theta_1, \ldots, \theta_d$ be a l.s.o.p. for $\field[\Delta]$.
Then 
$$\dim_\field (\field[\Delta]/(\Theta))_j
 = h_j(\Delta), \quad \mbox{ for all $0\leq j \leq d$.} 
$$ 
\end{theorem}

A generalization of Theorem \ref{CMh-numbers} for Buchsbaum 
complexes was found by Schenzel \cite{Sch}. To state it,
for a $(d-1)$-dimensional simplicial complex $\Delta$ 
and $0\leq j \leq d$, we define
\begin{equation} \label{h'}
h'_j(\Delta):= h_j + {d\choose j} 
\sum_{i=0}^{j-1} (-1)^{j-i-1} \beta_{i-1}(\Delta), \quad
\mbox{ where } \beta_{i-1} = \dim_\field \widetilde{H}_{i-1}(\Delta; \field).
\end{equation}
Note that if $\Delta$ is Cohen-Macaulay, then $h'_j(\Delta)=h_j(\Delta)$.

\begin{theorem} ({\bf Schenzel})   \label{Buchs-h'-numbers}
Let $\Delta$ be a $(d-1)$-dimensional Buchsbaum complex 
and let $\theta_1, \ldots, \theta_d$ be a l.s.o.p. for $\field[\Delta]$.
Then 
$$\dim_\field (\field[\Delta]/(\Theta))_j = h'_j(\Delta), 
\quad \mbox{ for all $0\leq j \leq d$.} 
$$ 
\end{theorem}

Using Hochster's formula and Schenzel's theorem, 
we can now derive a combinatorial-topological version of
our Theorem \ref{main-thm} for the case of Buchsbaum complexes 
as well as new lower bounds on their face numbers. 
This material concludes this section.

\begin{theorem}  \label{main2}
Let $\Delta$ be a $(d-1)$-dimensional Buchsbaum complex and 
let $\theta_1, \ldots, \theta_d$ be a l.s.o.p. for $\field[\Delta]$. 
Then for all $0\leq j \leq d$,
$$
\dim_\field \left(\soc \field[\Delta]/(\Theta)\right)_j 
\geq {d \choose j}\beta_{j-1}(\Delta).
$$
In particular, $h'_j(\Delta) \geq {d \choose j}\beta_{j-1}(\Delta)$, or, 
equivalently,
$h_j \geq {d \choose j}  \sum_{i=0}^{j} (-1)^{j-i} \beta_{i-1}(\Delta)$.
\end{theorem}
\proof For $0\leq j \leq d$, we have 
$$
\dim_\field (\soc \field[\Delta]/(\Theta))_j 
\geq {d \choose j} \dim H^j(\field[\Delta])_0 = {d \choose j}\beta_{j-1},
$$
where the first step follows from Theorem \ref{main-thm} 
and the second one from Hochster's formula. Since, 
$\soc \field[\Delta]/(\Theta) 
\subseteq \field[\Delta]/(\Theta)$,
Theorem \ref{Buchs-h'-numbers} completes the proof.
\endproof

\begin{corollary} \label{ker-ineq-cor}
Let $\Delta$ be a $(d-1)$-dimensional Buchsbaum complex,
let $\theta_1, \ldots, \theta_d$ be a l.s.o.p. for $\field[\Delta]$, and let 
$\omega\in \m_1$ be a linear form. Then for all $0< j \leq d$,
\begin{equation} \label{Ker-ineq}
\dim_\field \left(\ker \left(\field[\Delta]/(\Theta)_j
\stackrel{\cdot\omega}{\longrightarrow}
 \field[\Delta]/(\Theta)_{j+1}\right) \right) 
\geq {d \choose j}\beta_{j-1}(\Delta).
\end{equation}
\end{corollary}
\proof Use Theorem \ref{main2} and the fact that 
$\ker (\cdot \omega) \supseteq 
\soc \field[\Delta]/(\Theta)$.
\endproof

Corollary \ref{ker-ineq-cor} settles in the affirmative 
a part of \cite[Conjecture 36]{Kalai02} --- the conjecture that served as 
main motivation and starting point for this paper.

\section{Upper bounds on Buchsbaum complexes}  \label{Upper-section}
In this section we use Theorem \ref{main2} to 
derive new upper bounds on the face numbers of Buchsbaum simplicial complexes.
As an application, we prove K\"uhnel's conjecture on the Euler 
characteristic of even-dimensional manifolds.

Somewhat surprisingly, to describe the upper bounds 
on the $f$-numbers of simplicial complexes, one needs the notion of
a multicomplex.  A {\em multicomplex} $\M$ is a subset 
of monomials, say in variables
$x_1, \ldots, x_{n-d}$, that is closed under division, i.e.~if 
$\mu\in\M$ and $\nu | \mu$, then also $\nu\in\M$. For a multicomplex $\M$,
we denote by $\M_j$ the set of its elements of degree $j$, 
and by $F_j=F_j(\M)$ 
the cardinality of $\M_j$. We refer to 
$F(\M):=(F_0, F_1, \ldots)$ as the $F$-vector of $\M$.

The $F$-vectors of multicomplexes were completely characterized by 
Macaulay \cite{Mac} (see also \cite[Theorem II.2.2]{St96}).
Given two integers $l$ and $j$ there exists a unique expression of $l$ 
in the form
\begin{equation} \label{canonical}
l={n_j \choose j}+{n_{j-1} \choose j-1}+\cdots + {n_s \choose s}, 
\quad \mbox{ where } n_j>n_{j-1}>\cdots > n_s\geq s\geq 1.
\end{equation}
Define
$$
l^{<j>}:=
{n_j+1 \choose j+1}+{n_{j-1}+1 \choose j}+\cdots + {n_s+1 \choose s+1}.
$$

We say that $R$ is a  {\it standard graded $\field$-algebra} 
if it is  a $\Z$-graded $\field$-algebra with $R_i = 0$ for 
$i<0$, $R_0 \cong \field$ and is generated as an algebra by $R_1$ 
with $\dim_\field R_1 < \infty.$  
Equivalently, as a $\field$-algebra, $R \cong \field[x_1, \dots, x_m]/I$ 
for some homogeneous ideal $I.$ 
The Hilbert function of such an $R$ is the sequence 
$(\dim_\field R_0, \dim_\field R_1, \ldots).$

\begin{theorem} ({\bf Macaulay}) \label{Macaulay}
Let $\mathcal{F} = (F_0, F_1, \ldots)$ be a sequence of nonnegative integers. 
 The following are equivalent:
\begin{itemize}
\item
$\mathcal{F}$ is the $F$-vector of a 
multicomplex. 
\item
$F_0=1$ and 
$0\leq F_{j+1} \leq F_j^{<j>}$ for $j\geq 1$.
\item
$\mathcal{F}$ is the Hilbert function of a standard graded $\field$-algebra.  
\end{itemize}
\end{theorem}

Using Theorems \ref{CMh-numbers} and \ref{Macaulay}, 
Stanley \cite[Theorem 6]{St77} characterized all 
possible $h$-vectors of Cohen-Macaulay
simplicial complexes.

\begin{theorem} ({\bf Stanley}) \label{Stanley}
A vector $h=(h_0, h_1, \ldots, h_d)\in \Z^{d+1}$ is the $h$-vector
of a $(d-1)$-dimensional Cohen-Macaulay complex on $n$ vertices
if and only if $h_0=1$, $h_1=n-d$, and
$0\leq h_{j+1} \leq h_j^{<j>}$ for $1 \leq j\leq d-1$.
\end{theorem}

A generalization of the necessity portion of Theorem \ref{Stanley} for 
Buchsbaum complexes was given in \cite[Theorem 1.7]{N98}, 
where it was shown that if 
$\Delta$ is a $(d-1)$-dimensional Buchsbaum complex, then its $h'$-vector,
$(h'_0, h'_1, \ldots, h'_d)$, (as defined in (\ref{h'})) satisfies 
\begin{equation}      \label{old-bounds}
h'_{j+1}  \leq \left(h_j -{d-1 \choose j}\beta_{j-1}\right)^{<j>}, 
    \quad \mbox{ for } 1 \leq j\leq d-1.
\end{equation}
The first result of this section is to use Theorem \ref{main2} to 
strengthen the above inequalities.

\begin{theorem}  \label{new-bounds}
Let $\Delta$ be a $(d-1)$-dimensional Buchsbaum complex on $n$ vertices.
Then $h'_0=1$, $h'_1=n-d$, and
$$h'_{j+1}  \leq \left(h_j -{d \choose j}\beta_{j-1}\right)^{<j>}, 
   \quad \mbox{ for } 1 \leq j\leq d-1.
$$
\end{theorem}

\proof  Let $I = \soc \field[\Delta]/(\Theta).$ Since $I$ is killed by 
$\m$, any vector subspace of $I$ is an ideal of $\field[\Delta]/(\Theta).$  
In particular, $I_j$ is an ideal, so $(\field[\Delta]/(\Theta))/I_j$ 
is a standard graded $\field$-algebra. Let $(F_0, F_1, \dots, F_d, 0, \dots)$ 
be its Hilbert function. By Theorem \ref{Buchs-h'-numbers} and 
Theorem~\ref{main2}, 
$F_j \le h^\prime_j - \binom{d}{j} \beta_{j-1}$ 
and $F_{j+1} = h^\prime_{j+1}.$  Macaulay's theorem finishes the proof.
\endproof

The inequalities (\ref{old-bounds}) served in \cite{N98} as a key to
extending the Upper Bound Theorem for polytopes and spheres (UBT, for short)
to several classes of orientable homology manifolds
(among them, the class of all odd-dimensional homology manifolds and the class 
of all even-dimensional homology manifolds of Euler characteristic 2). 
This theorem,
originally proved by McMullen \cite{McMullen} for polytopes and 
later extended by Stanley \cite{St75} to homology spheres, 
asserts that among all
$d$-dimensional polytopes on $n$ vertices, the cyclic polytope, $C_d(n)$, 
has componentwise maximal $f$-vector.

A conjecture related to the UBT 
was proposed by K\"uhnel \cite[Conjecture B]{Kuh95}.
It asserts that if a simplicial complex $\Delta$ is a (combinatorial) 
$2k$-dimensional manifold (without boundary)
on $n$ vertices, then its Euler 
characteristic, 
$\chi(\Delta):=\sum_{j=0}^{2k} (-1)^j f_j=
1+\sum_{j=0}^{2k} (-1)^j \beta_j(\Delta)$, satisfies
\begin{equation} \label{Kuh-ineq}
(-1)^k (\chi(\Delta)-2) \leq \frac{{n-k-2 \choose k+1}}{{2k+1 \choose k}}.
\end{equation}
Moreover, equality happens if and only if $\Delta$ is {\em $(k+1)$-neighborly},
that is, every $k+1$ vertices of $\Delta$ form the vertex set of a 
face of $\Delta$.

While inequalities (\ref{old-bounds}) 
(together with  Klee's extension of the Dehn-Sommerville relations - 
Theorem \ref{klee} \cite{Klee})
were strong enough to imply the UBT
for several classes of homology manifolds, 
they were insufficient to completely prove the K\"uhnel conjecture, which 
was verified in \cite{N98} and \cite{N05} only for 
$2k$-dimensional orientable $\field$-homology manifolds with at least $4k+3$ or
at most $3k+3$ vertices. (Paper \cite{N98} treated the case of 
$\Char\field=0$, 
while \cite{N05} dealt with a field of an arbitrary characteristic.)
However, it was observed in \cite{N98} 
(see proof of Theorem 7.6 there) that 
if the inequalities of Theorem \ref{new-bounds} were true, they 
would imply K\"uhnel's conjecture for all $n$. Thus we now have

\begin{theorem} \label{Kuhnel-conj}
K\"uhnel's conjecture holds for all orientable 
$2k$-dimensional \field-homology manifolds.
In particular, K\"uhnel's conjecture holds 
for all  simplicial $2k$-manifolds.
\end{theorem}
\proof To prove the inequality, use Theorem \ref{new-bounds} 
together with \cite[Theorem 7.6]{N98}. The `in particular'-part 
follows from the fact that every simplicial manifold 
is orientable over a field of characteristic 2. 
Finally, the treatment of equality is given in \cite{N05} 
(see Remark at the end of Section 4 there.) \endproof

Using Theorem \ref{main2} and techniques developed in \cite{N05},
the inequalities of Theorem~\ref{new-bounds} 
can be significantly strengthened for the family of 
{\em centrally symmetric} Buchsbaum complexes, 
i.e.~complexes with a free $\Z/2\Z$-action. Combinatorially, these
inequalities can be described as follows.

\begin{theorem} \label{cs-new-bounds}
Let $\Delta$ be a $(d-1)$-dimensional 
centrally symmetric Buchsbaum complex with $n=2m$ vertices.
Then for every $1\leq j \leq d-1$, there exists a multicomplex 
$\M=\M(j)$ on $2m-d$ variables
$x_1, \ldots, x_{2m-d}$ all of whose elements are {\bf squarefree}
in the first $m$ variables and such that  
$$
F_{j+1}(\M)=h'_{j+1} (\Delta), \quad \mbox{ while } \quad 
F_{j}(\M)\leq h'_{j} (\Delta) - {d \choose j}\beta_{j-1} (\Delta).
$$
\end{theorem}

\proof Label the vertices of $\Delta$ so that for every $1\leq i \leq m$,
$x_i$ and $x_{m+i}$ are antipodal (i.e., $x_i, x_{m+i}$
form an orbit  under the given $\Z/2\Z$-action).
Consider $u\in \GL_n(\field)$ of the form
$$
u=\left[ \begin{array}{cc}
 I_m & I_m\\
  O  & Y^{-1}
         \end{array}
  \right].
\quad \mbox{ Equivalently, } \quad
u^{-1}=\left[ \begin{array}{cc}
 I_m & -Y\\
  O  & Y
         \end{array}
  \right].
$$
Here $I_m$ denotes the $m\times m$ identity matrix, $O$ stands for the 
$m\times m$ zero matrix, and $Y\in \GL_m(\field)$ satisfies the condition that
all of its $d \times d$-minors supported on the last $d$ columns of $Y$
are non-singular. Since $\field$ is infinite, such $Y$ exists.

Note that $u$ defines a graded automorphism 
of $S$ via $u(x_j)=\sum_{i=1}^n u_{ij}x_i$, and so
$uI_\Delta$ is a homogeneous ideal of $S$.
Let $I=uI_\Delta+(x_{n-d+1}, \ldots, x_n)$, 
let $\soc I=I:\m$ be the socle of $I$,
and let $J=I+(\soc I)_j$. 
Since for every element $y\in (\soc I)_j$, 
$\m \cdot y \subset I$,
it follows that $J$ is an ideal of $S$.
As no face of $\Delta$ contains two antipodal points,
the structure of $u^{-1}$ and \cite[Lemma III.2.4]{St96}  
imply that $x_{n-d+1}, \ldots, x_n$ is a l.s.o.p.~for $S/uI_\Delta$.
Hence, by Theorem \ref{Buchs-h'-numbers} and 
Theorem~\ref{main2}, 
\begin{equation}  \label{Hilb-S/J}
\dim_\field (S/J)_{j+1}=h'_{j+1} \quad \mbox{ and } \quad
\dim_\field (S/J)_{j} \le h^\prime_j - \binom{d}{j} \beta_{j-1}. 
\end{equation}

To construct a required multicomplex,
fix the reverse lexicographic order $\succ$ 
on the set of all monomials of $S=\field[x_1, \ldots, x_n]$ 
that refines the partial 
order by degree and satisfies $x_1\succ x_2 \succ \ldots\succ x_n$ 
(e.g.~$x_1^2\succ x_1x_2\succ x_2^2\succ x_1x_3 \succ x_2x_3 \succ x_3^2 
\succ \cdots$). Consider $\Init J$ --- the reverse lexicographic 
initial ideal of $J$ (see \cite[Section 15.2]{Eis}), 
and define $\M$ to be the collection of all monomials
that are not in $\Init J$. 
Since $\Init J$ is a monomial ideal that contains $x_{n-d+1}, \ldots, x_n$,
$\M$ is a multicomplex on $n-d$ variables. Moreover, 
$\M$ has ``correct'' $F$-numbers.
This follows from Eq.~(\ref{Hilb-S/J}) and the fact 
that $\M$ is a $\field$-basis 
for $S/J$ (see \cite[Thm.~15.3]{Eis}). Finally, the structure of $u$ and 
that $\{x_i, x_{i+m}\}$ is not a face
of $\Delta$ imply that
$\Init J\ni\Init u(x_ix_{i+m})=x_i^2$ for all $1\leq i \leq m$, and hence that
all elements of $\M$ are squarefree in the first $m$ variables.
\endproof

A complete characterization of $F$-vectors of multicomplexes that 
are squarefree in the first $m$ variables was worked out by
Clements and Lindsr{\"o}m \cite{ClLind}. Their theorem provides an
explicit sharp upper bound on $F_{j+1}$ of such a multicomplex in terms 
 of its $F_{j}$ and $j$. 
(Compare to Macaulay's theorem that characterizes $F$-vectors 
of multicomplexes without any restrictions on degrees.) 
Thus using Clements-Lindsr{\"o}m
theorem, one can restate Theorem \ref{cs-new-bounds} in purely numerical terms.

\begin{remark}{\rm
The same proof as in Theorem \ref{cs-new-bounds}
but with matrix $u$ chosen as in \cite[Theorem 3.3]{N05} allows to extend 
Theorem \ref{cs-new-bounds} to all Buchsbaum
simplcial complexes with a proper $\Z/p\Z$-action, where $p$ is a prime 
number, thus proving Conjecture 6.1 of \cite{N05}. So far we have been  
unable to settle Conjecture 6.2 of \cite{N05} --- an analog of K\"uhnel's
conjecture for manifolds with symmetry. The statement in \cite{N05} that 
\cite[Conj.~6.1]{N05} would imply \cite[Conj.~6.2]{N05} at least for all 
centrally symmetric manifolds is erroneous.}
\end{remark}

\section{Lower bounds} \label{lower_section}

The Dehn-Sommerville relations for simplicial
polytopes states that $h_i = h_{d-i}.$  
Klee proved an analogous formula for  semi-Eulerian complexes.  
A  pure complex is {\it semi-Eulerian} if the Euler characteristic of the 
link of every nonempty face is the same as the Euler characteristic of a 
sphere of the same dimension.  A prototypical example is an arbitrary
triangulation of a homology manifold without boundary.  

\begin{theorem}\cite{Klee}(Klee's formula) \label{klee}
Let $\Delta$ be a semi-Eulerian $(d-1)$-dimensional complex.  Then
$$h_{d-i} - h_i = (-1)^i {d \choose i} [\chi(\Delta) - \chi(S^{d-1})].$$
\end{theorem}

An immediate consequence of Klee's formula is that for semi-Eulerian complexes
 knowledge of the $g$-vector is sufficient to recover the $f$-vector.  
The {\it $g$-vector} of $\Delta$ is $(g_0, \dots, g_{\lfloor d/2 \rfloor})$, 
where $g_i = h_i - h_{i-1}.$  Of particular interest in this section is 
$g_2 = h_2 - h_1 = f_1 - d f_0 + \binom{d+1}{2}.$

In \cite{Kalai87} Kalai conjectured that if $\Delta$ is a triangulation of a 
closed manifold with $d \ge 4$, then 
$g_2 \ge \binom{d+1}{2} \beta_1(\Delta; \Q).$  
This bound is sharp for triangulations in $\H^d.$  A complex $\Delta$ is in 
$\H^d$ if it can be obtained from the boundary of the $d$-simplex by a 
sequence of the following three operations:

\begin{itemize}
\item
  Subdivide a facet with one new vertex in the interior of the facet.
\item
  Form the connected sum of $\Delta_1, \Delta_2 \in \H^d$ by identifying a 
pair of facets, one from each complex,  and then removing the interior of 
the identified facet.
  \item
  Form a handle by identifying a pair of facets in $\Delta \in \H^d$ 
and removing
 the interior of the identified facet in such a way that the resulting complex
 is still a simplicial complex.  Equivalently,  the distance in the 
$1$-skeleton between every pair of identified vertices must be at least three. 
\end{itemize}
If the only type of  operation used is the first one (subdividing a facet), then the resulting space is a {\it stacked sphere.} Another characterization of $\H^d$, due to Walkup in dimension three 
\cite{Wa} and Kalai in higher dimensions \cite{Kalai87}, 
is as those triangulations all of whose vertex links are stacked spheres.

Kalai's conjecture was verified for $\beta_1 = 1$ and for orientable manifolds 
when $d \ge 5$ and $\beta_2 = 0$ in \cite{Sw}. In the latter case, if 
$g_2 = \binom{d+1}{2} \beta_1(\Delta; \Q),$ then $\Delta \in \H^d.$  This last 
result was then used to determine all possible pairs $(f_0, f_1)$ for 
triangulations of spherical bundles over the circle \cite{CSS}.  

\begin{theorem} \label{g2}
Let $\Delta$ be a connected
triangulation of an orientable $\field$-homology $(d-1)$-dimensional 
manifold with $d \ge 4.$  Then 
\begin{equation}  \label{g2f}
g_2 \ge \binom{d+1}{2} \beta_1(\Delta; \field).
\end{equation} 
Furthermore, if $g_2 = \binom{d+1}{2} \beta_1$ and $d \ge 5,$ then 
$\Delta \in \H^d.$  
\end{theorem} 

\proof  First we consider the situation 
when the characteristic of $\field$ is zero.
By \cite[Lemma 5.1]{N98}, $h^\prime_{d-2} - h^\prime_2 = 
\binom{d}{2}(\beta_2 - \beta_1)$ and 
$h^\prime_{d-1} - h^\prime_1 = d (\beta_1 - \beta_0) = d \beta_1.$ 
For generic l.s.o.p. $\Theta$ and one-form $\omega,$ 
multiplication by  $\omega$ induces a surjection  
$\omega:(\field[\Delta]/( \Theta ))_{d-2} \to 
(\field[\Delta]/( \Theta ))_{d-1}$ \cite[Corollary 4.29]{Sw}.
Since the dimension of the socle of 
$(\field[\Delta]/( \Theta ))_{d-2}$ 
is at least $\binom{d}{d-2} \beta_{d-3},$

$$h^\prime_{d-2} - \binom{d}{d-2} \beta_{d-3} \ge h^\prime_{d-1}.$$

\noindent  Combining this with Poincar\'e duality
$$\begin{array}{ccl}
  h^\prime_2 + \binom{d}{2}(\beta_2 - \beta_1)  - 
       \binom{d}{2} \beta_2 & \ge & d \beta_1 + h^\prime_1 \\
  h^\prime_2 - h^\prime_1 &\ge &d \beta_1 + \binom{d}{2} \beta_1\\
  h_2 - h_1 & \ge & \binom{d+1}{2} \beta_1,
  \end{array}$$
where the last line follows from Schenzel's formula (see
Theorem \ref{Buchs-h'-numbers}).

Suppose  $g_2 = \binom{d+1}{2} \beta_1$ and $d \ge 5.$ 
The previous computation shows that 
for generic $\omega,$ the kernel of multiplication by $\omega$ in degree 
$d-2$ equals the socle of $\field[\Delta]/( \Theta )$ in degree 
$d-2.$  Consider the ideals generated by the variables $( x_i ).$ 
 By \cite{Sw}, $( x_i ) \subseteq 
\field[\Delta]/( \Theta )$ 
is isomorphic as an $S$-module to 
$\field[\lk i]/( \Theta^\prime )$ 
with a degree one shift for a suitably defined $\Theta^\prime.$  
Hence, if $\ker \omega \cap (( x_i ))_{d-2} \neq 0,$ 
then the socle of $(\field[\lk i]/( \Theta^\prime ))_{d-3}$ 
is also nonzero.  This is impossible since the link is a Gorenstein* complex. 
 As multiplication by a generic one-form from 
$(\field[\lk i]/( \Theta^\prime ))_{d-3}$ 
surjects onto $(\field[\lk i]/( \Theta^\prime ))_{d-2},\ 
h_{d-3}(\lk i) = h_{d-2} (\lk i).$ Equivalently, by the Dehn-Sommerville 
relations, $h_1( \lk i) = h_2 (\lk i).$ 
The lower bound theorem \cite[Theorem 1.1]{Kalai87} 
shows that each $\lk i$ must be a stacked sphere. 

What if the characteristic of $\field$ is not zero?  The only part of the 
above  which  needs to be changed is the proof that  for generic $\Theta$ 
and one-form $\omega,$ multiplication  induces a surjection  
$\omega:(\field[\Delta]/( \Theta ))_{d-2} \to 
(\field[\Delta]/( \Theta ))_{d-1}.$  The proof given in \cite{Sw} depends on
 \cite{Lee} and the generic rigidity of embeddings of polytopes in $\R^3.$  
Hence this approach is  only valid in characteristic zero.  However, Murai's 
recent preprint \cite[Corollary 3.5]{Mu}, combined with Whiteley's proof that 
two-dimensional spheres are strongly edge decomposable \cite{Wh} 
(see \cite{Nevo} for the definition of strongly edge decomposable), 
provide an alternative proof which is valid in nonzero characteristics. 
\endproof  

\begin{problem}
Suppose $\Delta$ is a $\field$-orientable $3$-dimensional manifold 
without boundary and $g_2 = 10 \beta_1.$  Is $\Delta \in \H^4?$
\end{problem}
When $\beta_1=1,$  the answer to this problem is known to be yes \cite{Wa}. 

Under certain conditions, Theorem \ref{main2}
can provide absolute 
lower bounds for $h^\prime$-vectors (and hence $f$-vectors) of Buchbaum 
complexes of a fixed homological type.  Suppose  
$\beta_{i-1} $ is the only nontrivial Betti number of $\Delta.$   
By Theorem \ref{main2}, $ h^\prime_i \ge \binom{d}{i} \beta_{i-1}.$  
Furthermore, assume that $\binom{d}{i} \beta_{i-1} = \binom{m}{i}$ 
for some $m.$ Macaulay's upper bound for Hilbert functions implies that for 
$j \le i, h^\prime_j \ge \binom{m-i+j}{j}.$   Thus,  if $\Delta$ satisfies all
 of these restrictions as equalities and $h^\prime_j = 0$ for $j>i,$  
then $\Delta$ has the minimum possible $h^\prime$-vector for a Buchbaum 
complex of this homological type.  Terai and Yoshida examined precisely 
this situation in \cite{TeYo}. 

\begin{theorem} {\rm \cite[Theorem 2.3]{TeYo}}  \label{TeYo}
  Let $\Delta$ be a $(d-1)$-dimensional Buchsbaum complex which is 
$i$-neighborly, but not $(i+1)$-neighborly.   
Set $\beta = \binom{n-d+i-1}{i}/\binom{d}{i}.$ 
Then the following are equivalent. 
\begin{itemize}
\item
  $\field[\Delta]$ has an $(i+1)$-linear resolution and $\beta_{i-1} = \beta.$
  \item
  The $h$-vector of $\Delta$ is 
  $$\left(1,n-d, \binom{n-d+1}{2},\dots,\binom{n-d+i-1}{i},
 -\binom{d}{i+1} \beta, \binom{d}{i+2} \beta, \dots, (-1)^{d-i} \beta\right).$$
  \item
  For every vertex $j,$ the link satisfies $h_m(\lk j) = 0$ 
if and only if $m > i-1.$
\end{itemize}
\end{theorem}

As the previous paragraph shows, any space satisfying the above conditions 
has the minimum $f$-vector among all Buchsbaum complexes with the specified 
$\beta_i.$  Terai and Yoshida also proved that Alexander duals of cyclic 
polytopes form an infinite family of examples of the above phenomenon.   

For examples with Betti numbers greater than one, let $\Delta$ be a 
$2k$-dimensional manifold which is also  $(k+1)$-neighborly.  
Now consider $\Delta$ with one vertex, say $n$,  and all of its incident 
faces removed and call this new complex $\Delta^\prime.$  
As $\Delta$ was $(k+1)$-neighborly, its only nonzero reduced  Betti numbers 
are $\beta_k$ and $\beta_{2k}.$  The Mayer-Vietoris sequence 
for $\Delta = \Delta^\prime \cup (n \ast \lk n)$ shows that the only nonzero 
reduced Betti number for $\Delta^\prime$ is $\beta_k.$  
Since $\Delta^\prime$ is a manifold with boundary it is Buchsbaum.  
The $h$-vector of the link of any vertex of $\Delta$ is given by 
$h_i = \binom{n-2k-2+i}{i},$ for $i \le k$ and $h_i = h_{2k-i}$ for 
$k < i \le 2k.$ Similarly, for each vertex $j <n$ the $h$-vector of 
$\overline{\Star} n \subset \lk j,$ the closed star of $n$ within the link of 
$j,$ is specified by the same equation for 
$i <k , h_i = h_{2k-i-1}$ for $k \le i \le 2k-1,$  and $h_{2k}=0.$   
Using the same reasoning as in \cite[Lemma 3]{Chari}, 

$$h_i(\lk_{\Delta^\prime} j) = h_i (\lk_\Delta j) - h_{2k-i}
(\overline{\Star} n \subset \lk j).$$

\noindent Hence $\Delta^\prime$ satisfies the third condition of 
Theorem \ref{TeYo}.

There are only a few known triangulations of $2k$-manifolds which are also 
$(k+1)$-neighborly.  For surfaces there are the $2$-neighborly triangulations 
 in \cite{JuRi} and \cite{Ri}.  Other examples include $\C P^2$ \cite{KuLa}, 
$K3$-surfaces \cite{CaKu}, $S^3 \times S^3$ \cite{Lutz}, and 
$\mathbb{H}P^2$ \cite{BrKu}, where $\mathbb{H}P^2$ is a manifold whose 
cohomology ring is isomorphic to the cohomology ring of the quaternionic 
projective plane.

\section{Buchsbaum simplicial posets}   \label{posets_section}
The goal of this section is to rework most of material of 
Section \ref{St-R-rings}, including Theorem~\ref{main2},
in the generality of Buchsbaum simplicial posets.
Simplicial posets (also sometimes referred
to in the literature as Boolean cell complexes or pseudo-simplicial
complexes) provide a certain generalization of simplicial complexes. 
We start by reviewing their definition and related notions 
as well as the corresponding algebraic background.

A {\em simplicial poset} is a (finite) poset $P$ that has a 
unique minimal element,  $\hat{0}$,
and such that for every $\tau\in P$, the interval 
$[\hat{0}, \tau]$ is a Boolean algebra \cite{St91}. In particular, $P$ is
graded, and the face poset of any simplicial complex is a simplicial poset. 
As with simplicial complexes, one can think of simplicial posets 
geometrically: 
it follows from results of \cite{Bj} that every simplicial poset 
$P$ is the face poset of a certain regular CW-complex, $|P|$,
all of whose cells are simplices and every two cells are attached 
along a possibly empty subcomplex of their boundaries 
(rather than just one face, as in the case of a simplicial complex). 
We call $|P|$ the {\em realization of $P$}, 
and refer to its elements as faces. 
It also follows from \cite{Bj} that $|P|$ has a 
well-defined barycentric subdivision which is the
simplicial complex isomorphic to the order
complex $\Delta(\overline{P})$ of the poset $\overline{P}=P- \{\hat{0}\}$.

As in the case of simplicial
complexes, we denote by $f_i=f_i(P)$ the number of 
$i$-dimensional faces of $|P|$ (equivalently, 
the number of rank $i+1$ elements of $P$), and by 
$f(P)=(f_{-1}, f_0, \ldots, f_{d-1})$ the {\em $f$-vector of $P$},
and we define the $h$-vector of $P$, 
$h(P)=(h_0, \ldots, h_d)$ according to Eq.~(\ref{h-vector}).
Here $d-1$ is the {\em dimension of $|P|$}, that is, the 
maximal dimension of faces of $|P|$. Equivalently, $d=\rk P$, the rank of $P$.
From now on we refer to $P$ and $|P|$ almost interchangeably.

As with simplicial complexes, 
we need a notion of a link: for an element $\tau$ of $P$, 
we define the {\em link of $\tau$ in $P$}, to be
$$\lk\tau=\lk_P(\tau):= \{\sigma\in P \ | \ \sigma\geq \tau\}.$$
It is easy to check that $\lk \tau$ is also a simplicial poset 
with its $\hat{0}$ element being $\tau$, and that if 
$F=\{\tau_0< \tau_1< \ldots<\tau_r=\tau\}$ is a saturated chain in 
$(\hat{0},\tau]$, then 
$\lk_{\Delta(\overline{P})}(F)\cong\Delta(\overline{\lk_P(\tau)})$.

Associated to a simplicial poset $P$  is an algebra 
$A_P$ \cite{St91}, defined as follows.
For each element $\tau$ of $P$, consider a variable $x_\tau$. 
Let $\tilde{S}$ be the polynomial ring 
$\field[x_\tau \ | \ \tau\in P]$. We assume that the set of atoms of $P$ 
(equivalently, the set of vertices of $|P|$) is $V(P)=[n]$, 
so that, $S=\field[x_1, \cdots, x_n]$ is a subring of $\tilde{S}$.
The {\em face ring of $P$}, $A_P$, is then $\tilde{S}/I_P$, where $I_P$ 
is the ideal of $\tilde{S}$ 
generated by the elements of the following form:
\begin{itemize}
\item $x_\tau x_\sigma$ for all pairs of elements $\tau, \sigma \in P$ that 
have no common upper bound in $P$.
\item $\x_\tau x_\sigma -\x_{\tau \wedge \sigma}\sum \x_\rho$  
     for pairs of  $\tau, \sigma$ 
     incomparable in $P$, where the sum is over the set of all 
{\em minimal upper bounds of $\tau$ and $\sigma$}. 
Note that if $\tau$ and $\sigma$ have an upper bound $\rho$, then
   $\tau \wedge \sigma$ is well-defined, as $\tau$ and $\sigma$ 
are elements of $[\hat{0}, \rho]$, a Boolean algebra.
\item $\hat{0}-1$. 
\end{itemize}
Defining $\deg x_\tau :=\rk \tau$ makes $A_P$ into a $\Z$-graded algebra. 
There is also a $\Z^n$-refinement of this grading on $A_P$ given by 
$\deg \tau:= \sum \{ e_i \ | \ i\in [n], i\leq\tau\}$.
Here $e_1, \ldots, e_n$ is the standard basis for $\Z^n$.

We cite from \cite{St91} a few basic properties of $A_P$:
\begin{itemize}
\item  $A_P$ is an algebra with straightening laws 
(this is  \cite[Lemma 3.4]{St91}). 
\item $A_P$ is integral over $S$ \cite[Lemma 3.9]{St91}. 
Since $A_P$ is also finitely-generated algebra over $S$,
it follows that $A_P$ is a (graded) Noetherian $S$-module.
\item The Krull dimension of $A_P$ is $\rk P=\dim P+1=:d$, and 
(as was the case for a simplicial complex)
the $\Z$-graded Hilbert series of $A_P$  is given by
$F(A_P, x)=(1-x)^{-d} \sum_{i=0}^d h_i(P) x^i$ (see \cite[Prop.~3.8]{St91}).
\end{itemize}

An analog of Hochster's formula for the local cohomology of $A_P$ 
(as a module over $S$)
was worked out by Duval in \cite[Theorem 5.9]{Du}.
\begin{theorem}{(\bf Duval)}  \label{Duval}
For a simplicial poset $P$ with $V(P)=[n]$, the $\Z^n$-graded 
Hilbert series of the local cohomology modules of $A_P$ as $S$-modules is
$$
F(H^i(A_P), \lambda)=\sum_{\tau\in P} \beta_{i-\rk(\tau)-1}(\lk \tau)
\prod_{j\in[n],\ j\leq\tau} 
\frac{\lambda_j^{-1}}{1-\mathbf{\lambda}_j^{-1}},
$$ 
where $\beta_i(\lk \tau)$ is the $i$th
reduced Betti number of the order complex $\Delta(\overline{\lk\tau})$ and 
$\lambda = (\lambda_1, \dots, \lambda_n).$
\end{theorem}

Call a simplicial poset $P$ a {\em Cohen-Macaulay poset} 
if its order complex, $\Delta(\overline{P})$, is a 
Cohen-Macaulay simplicial complex, as defined in Section \ref{St-R-rings}.
Similarly, call $P$ a {\em Buchsbaum poset} if $\Delta(\overline{P})$
is a Buchsbaum simplicial complex.
Stanley \cite[Cor.~3.7]{St91} showed that if $P$ is a Cohen-Macaulay 
simplicial poset, then its face ring, $A_P$, is Cohen-Macaulay 
as a module over itself or over $S$. 
Here we use Theorem \ref{Duval} to prove a similar result about
Buchsbaum posets.

\begin{proposition} \label{Buchs-poset}
If $P$ is a Buchsbaum simplicial poset, then the ring $A_P$ is Buchsbaum as 
an $S$-module.
\end{proposition}
 \proof Since $\Delta(\overline{P})$ is a Buchsbaum simplicial complex, say,
of dimension $d-1$, 
it follows from Theorem \ref{Reisner-Schenzel}, that for $i<d$, 
\begin{equation} \label{beta=0}
\beta_{i-\rk(\tau)-1}(\lk \tau)=0 \quad \mbox{ unless }\tau=\hat{0}.
\end{equation}
Thus, by Theorem \ref{Duval}, for $i<d$, 
$F(H^i(A_P), \mathbf{\lambda})=\beta_{i-1}(\Delta(\overline{P}))$ 
is a number rather 
than a series, and hence for $i<d$, $H^i(A_P)$ is concentrated in degree 0.
Therefore, $\m\cdot H^i(A_P)=0$. Also for $0\leq i < j <d$, the only 
integer degrees $p$ and $q$ for which
$(H^i(A_P))_p\neq 0$ and $(H^j(A_P))_q\neq 0$
are $p=q=0$. In particular, $0> i-j=(i+p)-(j+q)$, and so $(i+p)-(j+q)\neq 1$.
Proposition 3.10 on page 98 of \cite{StVo} then implies that $A_P$ is a 
Buchsbaum module. \endproof 

Stanley showed \cite[Section 3]{St91} that Theorem \ref{CMh-numbers}
 holds in the generality of Cohen-Macaulay simplicial posets, that is, 
if $P$ is a Cohen-Macaulay poset of rank $d$
and $\{\theta_1, \ldots, \theta_d\} \subset S$ is a l.s.o.p.~for $A_P$, then 
$\dim_\field (A_P/(\Theta)A_P)_j=h_j$ for all 
$0\leq j \leq d$. 

We next use Proposition \ref{Buchs-poset} to verify that
Schenzel's theorem, Theorem \ref{Buchs-h'-numbers}, also continues to
hold in the generality of Buchsbaum simplicial posets. Our proof mostly
mimics that of Schenzel and is included here only for completeness.

\begin{proposition}  \label{Buchs-h'-numebrs2}
Let $P$ be a rank $d$ Buchsbaum simplicial poset, let
$$
h'_j(P):=h_j(P)+{d \choose j}\sum_{i=0}^{j-1} (-1)^{j-i-1} 
 \beta_{i-1}(\Delta(\overline{P})) \qquad \mbox{ for } 0\leq j \leq d,
$$
and let
$\{\theta_1, \ldots, \theta_d\}\subset S$ be a l.s.o.p.~for $A_P$.
Then $\dim_\field(A_P/(\Theta)A_P)_j=h'_j$ for 
$0\leq j \leq d$.
\end{proposition}

\proof 
From the following exact sequence of graded $S$-modules:
$$
0 \longrightarrow (0 :_{A_P} \theta_1)(-1)
\longrightarrow A_P(-1) 
\stackrel{\cdot \theta_1}{\longrightarrow} A_P
\longrightarrow A_P/(\theta_1)A_P,
$$
we obtain an expression for the Hilbert series:
$$
(1-x)F(A_P, x) = F(A_P/(\theta_1)A_P, x)-x\cdot F((0 :_{A_P} \theta_1), x).
$$
 Iterating the above $d$ times yields
\begin{equation}  \label{Hilb-AP}
(1-x)^dF(A_P, x) = F(A_P/(\theta_1, \ldots, \theta_d)A_P, x) - 
                  \sum_{i=0}^{d-1} x (1-x)^i \cdot F(L_i, x),
\end{equation}
where $L_i:=((\theta_1, \ldots, \theta_{d-1-i})A_P : \theta_{d-i})/
             (\theta_1, \ldots, \theta_{d-1-i})A_P$.

Now, since $A_P$ is a Buchsbaum module (see Proposition \ref{Buchs-poset}),
we have
\begin{eqnarray*}
L_i &\cong& H^0(A_P/(\theta_1, \ldots, \theta_{d-1-i})A_P)  
              \qquad \; \mbox{ (by \cite[pp.~64-65]{StVo})} \\
    &\cong & \bigoplus_{l=0}^{d-1-i} {d-1-i \choose l}H^l(A_P)(-l)
             \quad \mbox{ (by \cite[Lemma II.4.14'(b)]{StVo})} \\
   &\cong& \bigoplus_{l=0}^{d-1-i} 
   \field^{{d-1-i \choose l}\beta_{l-1}(\Delta(\overline{P}))}(-l)
       \qquad \;\; \, \mbox{ (by Theorem \ref{Duval} and Eq.~(\ref{beta=0})}),
\end{eqnarray*}
and so
\begin{equation} \label{Hilb-Li}
F(L_i, x)=
\sum_{l=0}^{d-i-1} {d-i-1 \choose l} 
\beta_{l-1}(\Delta(\overline{P}))\cdot x^l.
\end{equation}
Plugging (\ref{Hilb-Li}) into (\ref{Hilb-AP}), and using that 
$F(A_P, x)=(1-x)^{-d}\sum_{i=0}^d h_i(P)x^i$ (see properties of $A_P$ listed
above in this section), completes the proof.
\endproof

We are now in a position to derive the following poset-generalization
of Theorem \ref{main2}.

\begin{theorem} \label{main-posets}
Let $P$ be a rank $d$ Buchsbaum simplicial poset and let
$\theta_1, \ldots, \theta_d$ be a l.s.o.p.~for $A_P$. Then for all
$0\leq j \leq d$,
$$ \dim_\field \left(\soc A_P/(\Theta)A_P\right)_j
 \geq {d \choose j} \beta_{j-1} (\Delta(\overline{P})).
$$
Hence, $h'_j(P)\geq                                                    
       {d \choose j} \beta_{j-1} (\Delta(\overline{P}))$, or, equivalently,
$h_j(P) \geq {d \choose j} \sum_{i=0}^j (-1)^{j-i} 
     \beta_{i-1}(\Delta(\overline{P})).$
\end{theorem}
\proof The proof is the same as in Theorem \ref{main2}, just use 
 Theorem \ref{Duval} instead of Theorem~\ref{Hochster} and  
Proposition \ref{Buchs-h'-numebrs2} instead of Theorem \ref{Buchs-h'-numbers}.
\endproof

\section{Examples, concluding remarks, and open problems} 
             \label{examples-section}
\subsection{Toward the $g$-conjecture}     

Perhaps the most important problem in the theory of $f$-vectors is the 
$g$-conjecture.  The most optimistic version states that if $\Delta$ is 
a $(d-1)$-dimensional $\field$-homology sphere and $\Theta$ is a l.s.o.p.~for 
$\field[\Delta],$ then for a generic one-form $\omega$ and $i \le d/2,$ 
multiplication 
$$
\omega^{d-2i}: \field[\Delta]/(\Theta)_i \to \field[\Delta]/(\Theta)_{d-i}
$$ 
is an isomorphism.  Kalai has suggested a far-reaching generalization of this 
to homology manifolds \cite{N98}.  
     
   Suppose $\Delta$ is a simplicial complex which 
is homeomorphic to a $\field$-orientable homology manifold.  Define
    $$h^{\prime\prime}_i = h^\prime_i - \binom{d}{i} \beta_{i-1}.$$
    \noindent As pointed out in \cite{N98}, 
$h^{\prime\prime}_{d-i} = h^{\prime\prime}_i$  for $1 \le i \le d-1. $ Let 
$$I = \bigoplus^{d-1}_{j=1} \soc (\field[\Delta]/(\Theta))_j.$$
Since $I$ is a vector subspace of the socle it is also an ideal of $\field[\Delta]/(\Theta).$ Now set
$\overline{\field[\Delta]} = 
(\field[\Delta]/\Theta)/I.$ 
By Theorem \ref{main2} the dimension of 
$\overline{\field[\Delta]}_i$ is at most $h^{\prime\prime}_i$  for $1 \le i \le d-1.$

    \begin{conjecture} \cite{N98}
    For generic $\omega \in \field[\Delta]_1$ and $1 \le i \le d/2,$
      \begin{itemize}
        \item
          $\dim_\field \overline{\field[\Delta]}_i = h^{\prime\prime}_i.$
      \item
      Multiplication 
     $\omega^{d-2i}: \overline{\field[\Delta]}_i 
     \to \overline{\field[\Delta]}_{d-i}$
     is an isomorphism.
     \end{itemize}  
    \end{conjecture}
    
 Consider the  special case of $\Delta$ a homology sphere. The first part of 
the above conjecture holds since $\field[\Delta]$ is Gorenstein*.  The 
second part is the  $g$-conjecture. This suggests the following  conjecture.   
  
 \begin{conjecture}
 Let $\Delta$ be a $(d-1)$-dimensional Buchsbaum complex over $\field.$  
Let $\S$ be given by Theorem \ref{main-thm}, with $M = \field[\Delta]$.  
Then $\dim_\field \S = \dim_\field \S_0=1$ if and only if $\Delta$ 
is an orientable 
$\field$-homology manifold without boundary. 
 \end{conjecture}
 
A closely related, but potentially weaker conjecture is the following.

\begin{conjecture}
Suppose $\Delta$ is a simplicial complex homeomorphic to a $(d-1)$-dimensional $\field$-homology manifold.  Then $\overline{\field[\Delta]}$ is a Gorenstein ring.
\end{conjecture}

 \subsection{How tight are the bounds?}   
    
Theorem \ref{main-posets} together with a complete characterization of the
$h$-numbers of Cohen-Macaulay simplicial posets,
\cite[Theorem 3.10]{St91}, naturally leads to the  following question.

\begin{question} \label{suff-cond}
Do the bounds
$h'_j \geq {d \choose j}\beta_{j-1}$ for $j=1,2, \cdots, d-1$
together with $h'_0=1$ and $h'_d=\beta_{d-1}$ generate the complete set
of sufficient conditions for the $h$-numbers of Buchsbaum simplicial posets
with prescribed Betti numbers?
\end{question}
We believe that the answer is yes, and hence that this set of conditions 
gives a complete characterization of the possible pairs $(h,\beta)$ for 
Buchsbaum simplicial posets. The following result provides  partial evidence.

\begin{proposition}
Let $b_1,\dots,b_{d-1}, h^\prime_1, \dots, h^\prime_{d-1}$ be nonnegative 
integers.  Assume $d \le 5$ or $b_2 = \dots = b_{d-3} = 0.$  Then there 
exists a Buchsbaum simplicial poset $P$ with $\beta_j(|P|) = b_j$ and 
$h^\prime_j(P) = h^\prime_j$ for all $1 \le j \le d-1$ if and only if 
$h^\prime_j \ge \binom{d}{j} b_{j-1}$ for all $1 \le j \le d-1.$
\end{proposition}

If $b_i=0$ for all $i \neq d-1$, then one can even find a shellable poset 
satisfying the conditions
of the proposition, see \cite[Theorem 3.10]{St91}. 
For the other combinations of $b_i$ satisfying the hypotheses,
the proposition is an immediate consequence of the next four lemmas.

\begin{lemma} \label{X(1,d)}
There exists a $(d-1)$-dimensional Buchsbaum simplicial poset $X=X(1,d)$
such that $\beta_1(X)=1,  \beta_i(X) = 0$ for $i \neq 1$ and
$$
h'_i(X)=\left\{\begin{array}{cc}
{d \choose i}, & \quad \mbox{ if $i=0$ or $i=2$} \\
0, & \quad \mbox{ otherwise}.
   \end{array}
\right.
$$
\end{lemma}

\proof One such $X$ is given by taking a stacked ball whose facets are
$$\{1,2,\ldots, d\}, \{2,3,\ldots, d+1\}, \ldots, \{d, d+1, \ldots, 2d-1\},$$
and identifying the codimension one face spanned by $1, 2, \ldots, d-1$
with the codimension one face spanned by  $d+1, d+2, \ldots, 2d-1$ 
(where vertex $i$ is identified with vertex $d+i$).  The realization of $X$, 
$|X|$, is a $(d-2)$-disk bundle over $\S^1$,  orientable or not depending on 
the parity of $d$.
Hence $X$ is a Buchsbaum simplicial poset satisfying 
$\beta_1(X)=1$ and $b_i = 0$ for $i \neq 1.$  
A straightforward computation shows that all $h'_i$ numbers of $X$ vanish
except for $h'_0$ and $h'_2$, which are equal to 1 and ${d \choose 2}$,
respectively.
\endproof

\begin{lemma} \label{X(d-2,d)}
There exists a $(d-1)$-dimensional Buchsbaum simplicial poset $X=X(d-2,d)$
such that $\beta_{d-2}(X)=1,  \beta(X)_i = 0$ for $i \neq d-2$ and
$$
h'_i(X)=\left\{\begin{array}{cc}
{d \choose i}, & \quad \mbox{ if $i=0$ or $i=d-1$} \\
0, & \quad \mbox{ otherwise}.
   \end{array}
\right.
$$
\end{lemma}

\proof
One possible construction for $X$ is as follows.  
The vertices of $X$ are $1,2,\dots,d.$  The $(d-3)$-skeleton of $X$ 
is the $(d-3)$-skeleton of the $(d-1)$-simplex.  For every subset 
of vertices of cardinality $d-1$ give $X$  two distinct codimension one faces.
Label these faces $A_1, A_2, \dots, A_d, B_1, B_2, \dots, B_d,$ where 
$A_i$ and $B_i$ are the two faces whose vertices do not contain $i.$  
Any potential facet of $X$ is described by choosing one of $A_i$ or $B_i$ for 
each $i$ as the boundary faces of the facet.   The facets of $X$ are the $d$ 
possible ways of choosing exactly one boundary face of type $B$ and 
the rest of type $A.$  

Since $X$ has the $(d-3)$-skeleton of the simplex and also contains the 
$(d-2)$-skeleton of the simplex, $\beta_i(X) = 0$ for $i < d-2.$  
It is easy to see that the kernel of the boundary map from the 
$(d-1)$-chains to the $(d-2)$-chains is zero, hence $\beta_{d-1}(X) = 0.$ 
 A check of the Euler characteristic of $X$ shows that $\beta_{d-2}=1.$ 
Now that the Betti numbers of $X$ are known, direct computation shows 
that $X$ has the required $h^\prime$ numbers.   
To see that $\tilde{H}_i (\lk \sigma)= 0$ for a face $\sigma$ and 
$i < d-|\sigma|-1,$ use the same argument, except that the kernel 
of the boundary map in dimension 
$(d-1-|\sigma|)$ is of dimension $|\sigma|-1.$  
\endproof

\begin{lemma} \label{glue}
Let $P_1$ and $P_2$ be two disjoint $(d-1)$-dimensional Buchsbaum 
simplicial posets. If $Q$ is obtained from $P_1$ and $P_2$ by identifying 
a facet of $P_1$ with that of $P_2$, then $Q$ is also a Buchsbaum poset. 
Moreover,
\begin{eqnarray} \label{b+b}
\beta_i(Q) &=&\beta_i(P_1) + \beta_i(P_2), \quad i=0,1,\ldots, d-1, \quad
\mbox{ and} \\
                \label{h+h}
h'_i(Q)&=&h'_i(P_1)+h'_i(P_2)  \quad i=1,2,\ldots, d.
\end{eqnarray}
\end{lemma}

\proof That the Betti numbers add when $P_1$ and $P_2$ are glued along a facet
is an easy application of the Mayer-Vietoris sequence and the fact that the
intersection of $P_1$ and $P_2$ is contractible. The same
Mayer-Vietoris sequence also shows that $Q$ is Buchsbaum. Since
$f_{i-1}(Q)=f_{i-1}(P_1)+ f_{i-1}(P_2) -{d \choose i}$, the defining relation
for the $h$-numbers implies that $h_i(Q)=h_i(P_1)+h_i(P_2)$ for $i\geq 1$,
which together with Eq.~(\ref{b+b}) yields (\ref{h+h}).
\endproof

\begin{lemma} Let $P$ be a $(d-1)$-dimensional Buchsbaum simplicial poset
and let $g'_1, \ldots, g'_d$ be nonnegative integers satisfying
$g'_i\geq h'_i(P)$ for all $i=1,\ldots,d$. Then there exists a
$(d-1)$-dimensional Buchsbaum simplicial poset $Q$ whose Betti numbers,
except possibly for $\beta_{d-1}$, coincide with those of $P$ and such that
$h'_i(Q)=g'_i$ for all $1\leq i \leq d$.
\end{lemma}

\proof By \cite[Theorem 3.10]{St91} there exists a shellable simplicial
poset $R$ such that $h_i(R)=g'_i-h'_i(P)$ for all $1\leq i \leq d$. 
Attaching $R$ to $P$ along a facet (as in the proof of Lemma~\ref{glue}) 
produces a required poset $Q$.
\endproof

In view of the last two lemmas, to answer Question \ref{suff-cond}
in the affirmative, it is enough to construct for every $d$ and $i\leq d-1$
a $(d-1)$-dimensional Buchsbaum simplicial poset $X=X(i,d)$ such that
$$
\beta_j(X)=\left\{\begin{array}{cc} 0, & \mbox{ if $j\neq i$} \\
                                    1,  & \mbox{ if $j= i$}
                 \end{array}
            \right.
\mbox{ and }
h'_j(X)=\left\{\begin{array}{cc} 0, & \mbox{ if $j\neq 0, i+1$} \\
                                 {d \choose j},  & \mbox{ otherwise.}
                 \end{array}
            \right.
$$
Lemmas \ref{X(1,d)} and \ref{X(d-2,d)} provide such a construction for 
$i=1$ and $i=d-2$ (and any $d$), 
$X(0,d)$ is the disjoint union of two $(d-1)$-simplices, while $X(d-1,d)$
can be obtained by gluing two $(d-1)$-simplices along their boundaries. 
A construction for $X(2,5)$ is also known.  A simplicial poset homeomorphic 
to $\C P^2$ with $h$-vector $(1,0,0,10,-5,2)$ is described in \cite{Ga}. 
Removing any facet (or more precisely, the open cell of a facet) is an example
 satisfying the requirements of $X(2,5).$

The problem of determining all possible $h$-vectors of Buchsbaum complexes (as opposed to posets) was previously considered by Terai \cite{Terai} and in dimension $2 \  (d=3)$ by Hanano \cite{Han}.  The linear inequalities established in \cite[Theorem 2.4]{Terai} also hold for Buchsbaum posets.  In fact, the stronger inequalities, $i~h_i + (d-i+1) h_{i-1} \ge 0, ~1 \le i \le d,$  hold for arbitrary pure posets whose vertex links have nonnegative $h$-vectors \cite[Proposition 2.3]{Sw-indbr}.  At this time we do not have enough examples to make a firm conjecture which determines all possible  $(h,\beta)$ pairs for Buchsbaum complexes.  Hence we finish with a final question.
 
 \begin{question}
 Are there other restrictions on pairs $(h,\beta)$ for $(d-1)$-dimensional Buchsbaum complexes other than those coming from Theorem \ref{main2}, 
 $$h^\prime_i \ge \binom{d}{i} \beta_{i-1},$$
 and Theorem \ref{new-bounds}
 $$h'_{i+1}  \leq \left(h_i -{d \choose i}\beta_{i-1}\right)^{<i>}?
$$
 \end{question}

\section*{Acknowledgements} Most of this work was done during the special 
semester (spring and summer 2007) at the Institute for Advance Studies in 
Jerusalem. We are grateful to IAS for hospitality and especially to Gil Kalai 
for organizing this semester. Our additional thanks go to Gil for 
encouragement during the period when we went from ``having a proof"-stage to 
``not having a proof"-stage and back quite a few times. We are also grateful 
to Eran Nevo who explained to us how to use Whiteley's paper \cite{Wh} 
to prove Theorem \ref{g2} in nonzero characteristics.

\end{document}